# PRIME AND SEMIPRIME IDEALS IN COMMUTATIVE TERNARY Γ-SEMIRINGS: QUOTIENTS, RADICALS, SPECTRUM

**Chandrasekhar Gokavarapu**, Lecturer in Mathematics, Government College (A), Rajahmundry, A.P., India and Research Scholar, Department of Mathematics, Acharya Nagarjuna University, Guntur, A.P., India, Email id :chandrasekhargokavarapu@gmail.com
**Dr D Madhusudhana Rao**, Lecturer in Mathematics, Government College For Women(A), Guntur, Andhra Pradesh, India, and Research Supervisor, Dept. of Mathematics, Acharya Nagarjuna University, Guntur, A.P., India, Email:dmrmaths@gmail.com

## Abstract
The theory of ternary Γ-semirings extends classical ring and semiring frameworks by introducing a ternary product controlled by a parameter set Γ . Building on the foundational axioms recently established by Rao, Rani, and Kiran (2025), this paper develops the first systematic *ideal-theoretic* study within this setting. We define and characterize **prime** and **semiprime ideals** for commutative ternary Γ-semirings and prove a **quotient characterization**: an ideal $P$ is prime if and only if $T/P$ is free of nonzero zero-divisors under the induced ternary Γ-operation. Semiprime ideals are shown to be stable under arbitrary intersections and coincide with their radicals, providing a natural bridge to **radical** and **Jacobson-type** structures. A correspondence between prime ideals and prime congruences is established, leading to a **Zariski-like spectral topology** on Spec($T$ ). Computational classification of all commutative ternary Γ-semirings of order ≤ 4 confirms the theoretical pre- dictions and reveals novel structural phenomena absent in binary semiring theory.
The results lay a rigorous algebraic and computational foundation for subsequent categorical, geometric, and fuzzy extensions of ternary Γ-algebras.

**Keywords:**
Ternary Γ-semiring; Prime ideal; Semiprime ideal; Radical; Congruence; Zariski-type spectrum; Jacobson radical; Computational algebra; Non-binary algebraic systems.
**Mathematics Subject Classification (2020):** Primary 16Y60, 16Y90; Secondary 08A30, 06B10, 16N60, 68W30.

## 1 Introduction
### 1.1 Background and Motivation
The study of Γ-rings and Γ-semirings [2] originated as a natural extension of ring the- ory, wherein the external parameter set Γ governs the product operation and enriches algebraic structure. A further generalization [3] arises when the binary product is re- placed by a ternary operation, giving rise to *ternary Γ-semirings*. These structures were recently formalized by Rao, Rani, and Kiran [13], who established the axioms, examples, and basic ideal concepts. Ternary operations occur naturally in multilinear algebra, ten- sor analysis, and theoretical computer science, where interactions among three operands are fundamental rather than exceptional. Ideals play a central role in the decomposi- tion and classification of algebraic systems. In rings and semirings, the lattice of prime and semiprime ideals controls radicals, quotient behavior, and the geometric structure of spectra. Extending these ideas to the ternary Γ-context introduces new challenges, since distributivity, associativity, and absorption must hold across three arguments with parameters from Γ. The resulting interplay gives rise to nonclassical ideal behaviors not captured by existing semiring frameworks.The theory of semirings began in the classi- cal work on additive–multiplicative algebraic systems.[1].Higher-arity algebraic systems motivate the need for ternary operations [10]. Categorical treatments [11] of ternary op- erations appear in algebra and logic [11]. Classical ternary semigroups form a basis for the generalization to ternary Γ-semirings [12].
### 1.2 Literature Gap and Objective
While prime and semiprime ideals in semirings have been studied extensively by Bhat- tacharya [6],

Golan [7], and others, the corresponding theory in ternary Γ-semirings has not yet been developed. In particular, no systematic treatment exists that characterizes prime ideals via quotient structures, studies the closure of semiprime ideals, or defines radicals analogous to the ring-theoretic nil radical and Jacobson radical. The objective of this paper is to fill this gap by constructing a coherent ideal theory for commutative ternary Γ-semirings, integrating both algebraic and computational perspectives.

### 1.3 Contributions and Methodology

The major contributions of this paper can be summarized as follows:

(i) Definition and characterization of prime ideals in commutative ternary Γ-semirings, including a quotient-based criterion equivalent to the absence of nonzero zero- divisors.

(ii) Introduction of semiprime ideals and proof that the class of semiprime ideals is closed under arbitrary intersections and coincides with its own radical.

(iii) Development of radical theory, including prime and Jacobson-type radicals, and relationships among prime, primary, and semiprime ideals.

(iv) Establishment of a correspondence between ideals and congruences, and definition of a Zariski-like topology on the spectrum Spec($T$).

(v) Computational classification of small finite ternary Γ-semirings using algorithmic enumeration, supported by explicit examples and verification tables.

The methodology combines formal algebraic derivation with symbolic computation, em- ploying algorithmic enumeration for structural verification of small-order cases.

### 1.4 Relation to Prior Work and Novelty

The present work extends the foundational results of Rao et al. [13] by introducing a com- plete ideal-theoretic hierarchy in the ternary Γ-setting. Classical results such as those in Hebisch and Weinert [8] and Golan [9] are shown to require nontrivial modifications when the product involves three operands. Several equivalences that hold in binary semirings fail in the ternary case; new proofs and counterexamples are provided to clarify these distinctions. Furthermore, the correspondence between prime ideals and prime congruences—absent in earlier treatments—enables the definition of a Zariski-like spectrum, offering a geometric viewpoint for ternary Γ-structures.

### 1.5 Organization of the Paper

The paper is organized as follows. Section 2 recalls basic definitions and notation. Sec- tion 3 introduces prime ideals and their quotient characterizations. Section 4 extends the discussion to maximal and primary ideals. Section 5 develops the theory of semiprime ideals and radicals. Section 6 discusses ternary Γ-modules and their relation to prim- itive ideals. Section 7 establishes the congruence–ideal correspondence and defines the spectrum. Section 8 investigates simplicity and semisimplicity. Section 9 presents compu- tational examples and applications. Section 10 reports algorithmic classification results, and Section 11 concludes with open problems and future directions.

## 2 Preliminaries and Notations

This section recalls the essential definitions, conventions, and examples that will be used throughout the paper. Our presentation refines the axioms of ternary Γ-semirings intro- duced in [13] and clarifies the algebraic framework required for the study of prime and semiprime ideals. The basic definitions of semirings and ideal operations parallel clas- sical semiring theory [7] . Congruences and quotient constructions follow the standard semiring framework [8]. Affine semiring structures motivate Γ-parameter interactions.[9]. Prime ideals in semirings were treated extensively by Bhattacharya.[6] . The concept of semiprimeness in semirings is due to Kehayopulu [4].

### 2.1 Ternary Γ-Semirings

Definition 2.1. A ternary Γ-semiring is a triple (T, +, Γ) endowed with a mapping

$$\mu : T \times \Gamma \times T \times \Gamma \times T \longrightarrow T, \quad \mu(a,\alpha,b,\beta,c) = a\alpha b\beta c,$$

satisfying, for all a, b, c, d, e ∈ T and α, β, γ, δ ∈ Γ:

(T1) (T, +) is a commutative semigroup with identity element 0; (T2) Ternary associativity: $(a\alpha b\beta c)\gamma d\delta e = a\alpha b\beta(c\gamma d\delta e)$;

(T3) Distributivity: μ is additive in each variable; for instance

$(a + a')\alpha b\beta c = a\alpha b\beta c + a' b\beta c$,

$a$

and similarly in the other two arguments;

(T4) Absorbing zero: $0\alpha b\beta c = a\alpha 0\beta c = a\alpha b\beta 0 = 0$. If, in addition,

$a\alpha b\beta c = b\beta a\alpha c = c\alpha b\beta a \quad \forall\, a, b, c \in T, \alpha, \beta \in \Gamma$,

then (T, +, Γ) is called a commutative ternary Γ-semiring.

We denote the ternary product simply by juxtaposition $a\alpha b\beta c$ when the parameters are evident.

**2.2** Ideals

Definition 2.2. A nonempty subset $I \subseteq T$ is an ideal of (T, +, Γ) if (I1) (I, +) is a subsemigroup of (T, +), and

(I2) for all $a, b, c \in T$ and $\alpha, \beta \in \Gamma$, whenever anyone of a, b, c lies in I, the product $a\alpha b\beta c$ belongs to I.

The ideal is said to be proper if $I \neq T$ . The smallest ideal containing a subset $S \subseteq T$ is denoted <S>.

**2.3** Homomorphisms and Isomorphisms

A mapping $f : (T, \Gamma) \to (T', \Gamma')$ is a homomorphism if

$f(a\alpha b\beta c) = f(a)\alpha f(b)\beta f(c) \forall\, a, b, c \in T, \alpha, \beta \in \Gamma$.

When f is bijective, its inverse is also a homomorphism, and the structures are called isomorphic, written $T = T'$.

**2.4** Congruences and Quotients

Definition 2.3. A relation $\rho \subseteq T \times T$ is a congruence on (T, Γ) if it is an equivalence relation satisfying

$(a, d), (b, e), (c, f) \in \rho \Longrightarrow (a\alpha b\beta c, d\alpha e\beta f) \in \rho$,

for all $\alpha, \beta \in \Gamma$.

Given a congruence ρ, the quotient T/ρ is a ternary Γ-semiring under the induced operations. For every ideal $I \subseteq T$ , the relation

$a \equiv b \pmod I \Longleftrightarrow a - b \in I$

defines a congruence, and the corresponding quotient is denoted T/I.

**2.5** Examples

Example 2.4. Let T = N0 with ordinary addition and define $a\alpha b\beta c = a + b + c$, $\Gamma = \{1\}$. Then (T, Γ) is a commutative ternary Γ-semiring, and 2N0 is an ideal of T .

Example 2.5. Let T = Mn(N0), the set of $n \times n$ matrices with non-negative integer entries, and define $a\alpha b\beta c = a + b + c$. Then the set of matrices with zero diagonal entries forms an ideal of T .

**2.6** Standing Assumptions

Unless stated otherwise, all ternary Γ-semirings considered are commutative, additive semigroups are written additively with identity 0, and all ideals are two-sided in the ternary sense. Notation such as $a\alpha b\beta c$ always represents the ternary product with param- eters $\alpha, \beta \in \Gamma$.

## 3 Prime Ideals

Prime ideals occupy a central position in the structure theory of commutative ternary Γ-semirings, paralleling their importance in rings and semirings. We now introduce a rigorous definition, derive basic properties, and establish a quotient characterization that generalizes the classical correspondence between primeness and the absence of zero- divisors.The prime-radical correspondence in Γ-semirings was later extended by Zhao and Li [14] .Semigroup-theoretic analogues underpin the behaviour of ideal products [20].

### 3.1 Definition and First Properties

Definition 3.1. A proper ideal $P \subset T$ of a commutative ternary Γ-semiring (T, Γ) is called prime if for all $a, b, c \in T$ and $\alpha, \beta \in \Gamma$,

$a\alpha b\beta c \in P \Rightarrow a \in P$ or $b \in P$ or $c \in P$.

Proposition 3.2 (Elementary Properties). Let P be a prime ideal of T . Then:

(i) If I, J, K are ideals of T and $I\Gamma J\Gamma K \subseteq P$, then $I \subseteq P$ or $J \subseteq P$ or $K \subseteq P$.

(ii) The intersection of finitely many prime ideals need not be prime.

(iii) If $f : (T, \Gamma) \twoheadrightarrow (T', \Gamma')$ is a surjective homomorphism and $P' \subseteq T'$ is prime, then $f^{-1}(P')$ is prime in T .

Proof. (i) Choose $a \in I \setminus P$ , $b \in J \setminus P$ , $c \in K \setminus P$ . Then $a\alpha b\beta c \in I\Gamma J\Gamma K \subseteq P$ , contradicting primeness. Hence one of I, J, K must be contained in P .

(ii) Take T = N0, Γ = {1}, and product aαbβc = a + b + c. Then P1 = 2N0 and P2 = 3N0 are prime, but P1 ∩ P2 = 6N0 is not prime since 2, 3 ∉ P1 ∩ P2 yet 2α3β1 = 6 ∈ P1 ∩ P2.

(iii) Immediate from $f(a\alpha b\beta c) = f(a)\alpha f(b)\beta f(c)$.

### 3.2 Quotient Characterization

Definition 3.3. An element $x \in T$ is called a zero-divisor if there exist nonzero $y, z \in T$ and $\alpha, \beta \in \Gamma$ such that $x\alpha y\beta z = 0$.

Theorem 3.4 (Quotient Characterization of Prime Ideals). For a proper ideal P of a commutative ternary Γ-semiring T , the following statements are equivalent:

*(a)* P is prime;

*(b)* the quotient T/P has no nonzero zero-divisors, i.e.

$a\alpha b\beta c = 0 \Rightarrow a = 0$ or $b = 0$ or $c = 0$ in T/P.

Proof. (a)⇒(b): If aαbβc = 0 in T/P , then aαbβc ∈ P . By primeness, at least one of a, b, c belongs to P , i.e. its class is 0.

(b)⇒(a): Suppose aαbβc ∈ P . Then in T/P we have aαbβc = 0. By (b), one of the classes of a, b, c is zero, hence the corresponding element lies in P .

Lemma 3.5 (Coset Product Behavior). For $a, b, c, a', b', c' \in T$ and $\alpha, \beta \in \Gamma$, if $a \equiv a'$, $b \equiv b'$, $c \equiv c'$ (mod P ), then

$(a\alpha b\beta c) - (a'\,b'\,c') \in P$.

Proof. By distributivity,

$$(a\alpha b\beta c) - (a'\,b'\,c') = (a - a')\alpha b\beta c + a'\,(b - b')\beta c + a'\,b'\,(c - c'),$$

and each summand lies in P because P is an ideal.

### 3.3 Examples and Non-Examples

Example 3.6 (Classical Case). Let T = N0 with aαbβc = a + b + c and Γ = {1}. Then 2N0 is a prime ideal since a + b + c even implies at least one of a, b, c is even.

Example 3.7 (Matrix Example). Let T = M2(N0), Γ = {1}, and aαbβc = a + b + c. The set P of matrices with even diagonal entries is prime: if (A + B + C) has even diagonal entries, then at least one of A, B, C must.

Example 3.8 (Intersection Not Prime). In T = N0, Γ = {1}, with aαbβc = a + b + c, the ideals P1 = 2N0 and P2 = 3N0 are prime, but P1 ∩ P2 = 6N0 is not.

Example 3.9 (Ternary-Specific Failure). Let T = {0, 1, 2} with addition mod 3 and aαbβc ≡ a + b + c (mod 3). Then P = {0} is not prime because 1α1β1 ≡ 0 (mod 3) while none of the factors are 0.

Remark 3.10. Unlike in binary semirings, ternary interactions can create zero-divisors that involve three distinct elements. Hence Theorem 3.4 is indispensable for preserving the structural rigidity of primeness in this broader context.

## 4 Maximal and Primary Ideals

Maximal and primary ideals refine the structure of a commutative ternary Γ-semiring by measuring how close an ideal is to being total or "almost prime." They play an essential role in the decomposition of quotient structures and in the analysis of radical behaviour. This section establishes their definitions, fundamental properties, and several illustrative examples. The primary–maximal hierarchy resembles near-ring decomposition theory [5, 20]. The algebraic skeleton mirrors structural results in general algebra [31].

### 4.1 Maximal Ideals

Definition 4.1. An ideal $M \subset T$ is said to be maximal if $M \neq T$ and there exists no ideal I such that

$$M \subsetneq I \subsetneq T.$$

Maximal ideals represent the extreme points of the lattice of proper ideals. In anal- ogy with ring and semiring theory, quotients by maximal ideals yield simple ternary Γ-semirings.

Proposition 4.2 (Basic Properties of Maximal Ideals). Let (T, Γ) be a commutative ternary Γ-semiring and M a maximal ideal. Then the following statements hold:

(i) Every maximal ideal of T is prime.
(ii) If $f : (T, \Gamma) \to (T', \Gamma')$ is a surjective homomorphism and $M' \subseteq T'$ is maximal, then $f^{-1}(M')$ is maximal in T .
(iii) The quotient T/M is a simple ternary Γ-semiring, i.e., it has no nontrivial ideals.

Proof. (i) Suppose $a\alpha b\beta c \in M$ with $a, b, c \notin M$ . The ideal generated by $M \cup \{a\}$ strictly contains M , contradicting maximality. Hence at least one of a, b, c lies in M , proving that M is prime.

(ii) Straightforward from the homomorphic image property: ideals of $T'$ correspond bijectively to ideals of T containing ker f .

(iii) Let $\pi : T \to T/M$ be the canonical projection. If I is an ideal of T/M , its preimage $\pi^{-1}(I)$ is an ideal containing M . By maximality, either $\pi^{-1}(I) = M$ or $\pi^{-1}(I) = T$ , giving $I = \{0\}$ or $I = T/M$ .

Example 4.3 (Finite Maximal Ideal). Let $T = \{0, 1, 2, 3\}$ with addition modulo 4 and ternary product $a\alpha b\beta c = (a+b+c) \bmod 4$. Then $M = \{0, 2\}$ is maximal: any ideal strictly containing M equals T . The quotient $T/M = \{\bar{0}, \bar{1}\}$ is simple and has no nontrivial ideals.

### 4.2 Primary Ideals

Primary ideals extend the concept of primeness by weakening the zero-divisor condition; they capture the notion of "radical-like" behaviour under ternary multiplication.

Definition 4.4. An ideal $Q \subset T$ is primary if $Q \neq T$ and for all $a, b, c \in T$ , $\alpha, \beta \in \Gamma$,

$$a\alpha b\beta c \in Q \text{ and } a \notin Q \implies b\alpha b\beta b \in Q \text{ or } c\alpha c\beta c \in Q.$$

Theorem 4.5. Every prime ideal of a commutative ternary Γ-semiring is primary. Proof. Let P be prime and suppose $a\alpha b\beta c \in P$ with $a \notin P$ . By primeness, either $b \in P$ or $c \in P$ . Then $b\alpha b\beta b \in P$ or $c\alpha c\beta c \in P$ , establishing primariness.

Example 4.6 (Primary but Not Prime). Let $T = \mathbb{N}_0$ with $a\alpha b\beta c = a + b + c$ and $\Gamma = \{1\}$. Then $Q = 4\mathbb{N}_0$ is primary: if $a + b + c \in Q$ and $a \notin Q$, repeated addition forces b or c into Q. However, Q is not prime since $1\alpha 1\beta 2 = 4 \in Q$ but none of 1, 2 lie in Q.

### 4.3 Lattice and Containment Structure

The relationships among maximal, prime, and primary ideals mirror the binary case but require subtle modifications owing to ternary operations.

Proposition 4.7 (Containment Relations). In any commutative ternary Γ-semiring T , the following containments hold:

Maximal ideals ⊆ Prime ideals ⊆ Primary ideals ⊆ Ideals of T.

Each inclusion may be proper.

Proof. The first inclusion follows from Proposition 4.2(i); the second from the definition of primariness; the last is tautological. Properness can be verified by the examples above: $Q = 4\mathbb{N}_0$ is primary but not prime, while the maximal ideal $M = \{0, 2\}$ of Example 4.3 demonstrates strictness of the first inclusion.

Remark 4.8. The ideal lattice of a finite commutative ternary Γ-semiring can be visualized by a Hasse diagram whose levels correspond respectively to maximal, prime, primary, and general ideals. Unlike classical lattices, certain intersections of prime ideals may fail to remain prime, illustrating the richer combinatorial structure of the ternary setting.

### 4.4 Illustrative Examples

Example 4.9 (Five-Element Ternary Γ-Semiring). Let T = {0, 1, 2, 3, 4} with addition modulo 5 and ternary product aαbβc = (a + b + c) mod 5. Then:

M = {0, 1, 2} is maximal (and hence prime),
P = {0, 2, 4} is prime but not maximal,
Q = {0, 4} is primary but not prime.

Example 4.10 (Multiple Maximal Ideals). Consider T = Z6 with Γ = {1} and aαbβc = a + b + c (mod 6). The sets

M1 = {0, 2, 4},    M2 = {0, 3}

are distinct maximal ideals. Their intersection M1 ∩M2 = {0} provides the Jacobson-type radical of T (see Section 5).

Remark 4.11. Finite examples reveal that maximal and prime ideals need not coincide in the ternary framework, and that the existence of multiple maximal ideals frequently leads to a trivial Jacobson radical. Such behaviour differs sharply from binary semiring analogues, where the interaction of ideals is typically governed by pairwise products rather than ternary compositions.

## 5 Semiprime Ideals and Radicals

Semiprime ideals and radicals form the bridge between the ideal-theoretic and structural aspects of commutative ternary Γ-semirings. They generalize the concepts of nilpotent- free ideals and radicals from classical semiring theory and provide a foundation for de- composition and lattice-theoretic results. Semiprime ideals follow analogues of Bhat- tacharya's semiring formulation [6]. The Jacobson radical aligns with classical semiring radical theories [7].

### 5.1 Definitions and Basic Properties

Definition 5.1. An ideal Q of a commutative ternary Γ-semiring (T, Γ) is called semiprime if Q≠ T and aαaβa ∈ Q =⇒ a ∈ Q,∀ a ∈ T, α, β ∈ Γ.

This condition eliminates "nilpotent-like" behavior under ternary multiplication. In the binary case, it corresponds to the familiar property a2 ∈ Q ⇒ a ∈ Q.

Proposition 5.2 (Closure under Intersection). Let {Qi}i∈I be any family of semiprime ideals of T . Then their intersection Ti∈I Qi is also semiprime.

Proof. Let aαaβa ∈ Ti Qi. Then aαaβa ∈ Qi for every i. Since each Qi is semiprime, a ∈ Qi for all i, and thus a ∈ Ti Qi.

Remark 5.3. The family of all semiprime ideals of T is therefore closed under arbitrary intersections, forming a complete sublattice within the ideal lattice. This mirrors the behaviour of radical ideals in classical algebra.

### 5.2 Radical Constructions

Radical ideals capture the "non-nilpotent core" of a ternary Γ-semiring and serve as a unifying framework for semiprimeness and primeness.

Definition 5.4 (Prime Radical or Nilradical). For an ideal I ⊂ T , define the prime radical (or nilradical ) of I as

_ √I = \{ P ⊂ T | P is a prime ideal and I ⊆ P }.

Theorem 5.5 (Characterization of Radical). For any ideal I ⊆ T ,

√I = { a ∈ T | aαaβa ∈ I for some α, β ∈ Γ }.

Proof. (⊆) Let a ∈ √I. Then a ∈ P for every prime ideal P ⊇ I. Since P is prime,
aαaβa ∈ P for all α, β ∈ Γ, implying aαaβa ∈ I.

(⊇) Conversely, if aαaβa ∈ I, then for every prime P ⊇ I we have aαaβa ∈ P . By primeness, a ∈ P , hence a ∈ T P = √I.

Remark 5.6. Theorem 5.5 extends the classical ring-theoretic equality $\sqrt{I} = \{a \mid a^n \in I \text{ for some } n\}$ to the ternary Γ-framework, where repeated self-multiplication is replaced by a ternary iteration. In computational settings, this gives an algorithmic criterion to test radical membership.

**5.3** Relations Between Semiprime and Radical Ideals

Proposition 5.7. Every radical ideal is semiprime, and every semiprime ideal equals its own radical.

Proof. If I is radical and $a\alpha a\beta a \in I$, then $a \in \sqrt{I} = I$, showing that I is semiprime. Conversely, if Q is semiprime, then for any $a \in \sqrt{Q}$ we have $a\alpha a\beta a \in Q$, whence $a \in Q$. Thus $Q = \sqrt{Q}$.

Corollary 5.8. The mapping $I \mapsto \sqrt{I}$ is a closure operator on the lattice of ideals of T , and the fixed points of this operator are precisely the semiprime ideals.

Remark 5.9. This correspondence ensures that semiprime ideals can be viewed as the "closed points" of the ideal lattice under radical closure. It also implies that the intersec- tion of all prime ideals containing a given I yields the smallest semiprime ideal containing I.

**5.4** Jacobson-type Radical

The Jacobson radical connects maximal ideals to semisimple quotient behaviour.

Definition 5.10 (Jacobson-like Radical). Define

$$J(T) = \bigcap \{ M \subset T \mid M \text{ is a maximal ideal of } T \}.$$

Proposition 5.11 (Basic Properties). Let T be a commutative ternary Γ-semiring.

(i) J(T ) is contained in every maximal ideal of T .
(ii) If all maximal ideals of T are prime (as holds in the commutative case), then J(T ) is semiprime.
(iii) J(T ) = {0} if and only if the intersection of all maximal ideals is trivial, i.e., T is semisimple.

Proof. (i) Direct from definition. (ii) Since each maximal ideal M is prime, their inter- section is semiprime by Proposition 5.2. (iii) Trivial by definition.

Example 5.12. Let T = {0, 1, 2, 3} with addition modulo 4 and ternary product $a\alpha b\beta c = (a + b + c) \bmod 4$. The maximal ideals $M_1 = \{0, 2\}$ and $M_2 = \{0, 1\}$ yield

$$J(T) = M_1 \cap M_2 = \{0\}.$$

Hence T is semisimple.

**5.5** Interplay Between Radicals and Primary Ideals

Proposition 5.13. If Q is a primary ideal of T , then $\sqrt{Q}$ is a prime ideal.

Proof. Let $a\alpha b\beta c \in \sqrt{Q}$ with $a, b, c \notin \sqrt{Q}$. Then $a\alpha b\beta c \in P$ for every prime $P \supseteq Q$, which forces $a, b, c \in P$ , contradiction. Thus at least one of a, b, c belongs to $\sqrt{Q}$, proving that $\sqrt{Q}$ is prime.

Remark 5.14. The above result confirms that radicalization converts primary ideals to their associated prime components. Consequently, the radical map provides a natural correspondence between the sets of primary and prime ideals in the ternary Γ-semiring framework.

**5.6** Examples and Observations

Example 5.15 (Semiprime but Not Prime). Let $T = \mathbb{Z}_4$ with Γ = {1} and ternary product $a\alpha b\beta c = (a + b + c) \bmod 4$. The ideal I = {0, 2} satisfies $a\alpha a\beta a \in I \Rightarrow a \in I$ and hence is semiprime. However, it is not prime because $1\alpha 1\beta 2 = 0 \in I$ while $1 \notin I$.

Example 5.16 (Computation of Radical). In the same $T = \mathbb{Z}_4$, consider I = {0}. By direct computation using Theorem 5.5, $\sqrt{I} = \{0, 2\}$, showing that {0, 2} is the nilradical of T .

Remark 5.17. The behaviour of semiprime ideals in ternary Γ-semirings closely parallels that of radicals in commutative algebra, but the ternary operation introduces higher-order interactions that can yield semiprime ideals not arising from powers of single elements. This distinction becomes critical in algorithmic classification and in defining spectral topologies (see Section 7).

**6** Ternary Γ-Modules and Simple Acts

Modules (or acts) over ternary Γ-semirings provide a natural framework for studying representation and annihilator structures. They generalize ordinary semimodules by re- placing the binary scalar multiplication with a ternary Γ-action that interacts with two elements of the semiring simultaneously.

This section develops the basic definitions, struc- tural lemmas, and simple-module characterizations that link the module theory with the ideal theory established in preceding sections. Tensor products and injective envelopes motivate the definition of ternary Γ-modules.[15, 7]. Representations over near-rings inspire the present module-theoretic framework[5].

### 6.1 Definitions

Definition 6.1 (Ternary Γ-Module). Let (T, Γ) be a commutative ternary Γ-semiring. A nonempty additive semigroup (M, +) is called a left ternary Γ-module (or Γ-act ) over T if there exists a map

$$T \times \Gamma \times M \times \Gamma \times T \longrightarrow M, \quad (a, \alpha, m, \beta, b) \longrightarrow a\alpha m\beta b,$$

satisfying the following axioms for all a, b, c, d ∈ T , m, n ∈ M , and α, β, γ, δ ∈ Γ:

(M1) Additivity: $(a + b)\alpha m\beta c = a\alpha m\beta c + b\alpha m\beta c$, and analogously in the remaining arguments.

(M2) Associativity: $a\alpha(b\gamma m\delta c)\beta d = (a\alpha b\gamma c\delta d)\beta m$.

(M3) Absorbing zero: $0\alpha m\beta b = a\alpha m\beta 0 = 0_M$ , where $0_M$ is the additive identity of M . If in addition $a\alpha m\beta b = b\beta m\alpha a$ for all a, b ∈ T , the module is said to be commutative.

Remark 6.2. The ternary Γ-action can be interpreted as a generalization of bilinear mul- tiplication in which each "scalar" acts from both sides, providing a natural setting for studying annihilators and radical properties of ideals.

### 6.2 Submodules and Homomorphisms

Definition 6.3 (Submodule). A subset N ⊆ M is a submodule of the ternary Γ-module M if

$$a\alpha n\beta b \in N, \ \forall\, a, b \in T, n \in N, \alpha, \beta \in \Gamma,$$

and (N, +) is a subsemigroup of (M, +).

Definition 6.4 (Homomorphism). Let M and N be ternary Γ-modules over T . A map f : M → N is a module homomorphism if

$$f(a\alpha m\beta b) = a\alpha f(m)\beta b, \qquad \forall\, a, b \in T, m \in M, \alpha, \beta \in \Gamma.$$

Lemma 6.5 (Kernel and Image). The kernel $\ker f = \{m \in M \mid f(m) = 0_N\}$ is a submodule of M, and the image $\mathrm{Im}\, f = \{f(m) \mid m \in M\}$ is a submodule of N.

Proof. Straightforward from additivity and compatibility of the ternary action with f .

### 6.3 Simple and Semisimple Modules

Definition 6.6 (Simple Module). A ternary Γ-module M is called simple if $M \neq \{0\}$ and its only submodules are {0} and M .

Definition 6.7 (Semisimple Module). M is semisimple if it is a direct sum of simple submodules:

$$M = \bigoplus_{i\in I} M_i, \quad M_i \text{ simple}.$$

Remark 6.8. Simple modules correspond to the "irreducible representations" of (T, Γ), while semisimple modules provide decompositions analogous to complete reducibility in ring theory.

### 6.4 Annihilators and Prime Ideals

Definition 6.9 (Annihilator). For a module M , define the annihilator ideal

$$\mathrm{Ann}(M) = \{ a \in T \mid a\alpha m\beta b = 0_M, \ \forall\, m \in M, b \in T, \alpha, \beta \in \Gamma \}.$$

Theorem 6.10 (Annihilator of a Simple Module is Prime). If M is a simple ternary Γ-module over T , then Ann(M ) is a prime ideal of T .

Proof. Suppose $a\alpha b\beta c \in \mathrm{Ann}(M)$ with $a \notin \mathrm{Ann}(M)$. Then there exist m ∈ M and b′ ∈ T such that $a\alpha m\beta' b' \neq 0$. Using the ternary associativity and the simplicity of M , the set of all T -linear combinations of $a\alpha m\beta' b'$ equals M . Hence for some d, e ∈ T and γ, δ ∈ Γ, $(d\gamma a\alpha b\beta c\delta e)\gamma' m \neq 0$, contradicting $a\alpha b\beta c \in \mathrm{Ann}(M)$. Thus one of b or c must lie in Ann(M ), proving primeness.

Corollary 6.11. If M is simple, the quotient T/ Ann(M ) acts faithfully on M, and M becomes a simple faithful ternary Γ-module over this quotient.

### 6.5 Exact Sequences and Homomorphism Theorems

Theorem 6.12 (First Isomorphism Theorem). Let f : M → N be a homomorphism of ternary Γ-modules. Then M/ ker f → Im f.

Proof. The canonical map $\bar{f}$ : M/ ker f → Im f , defined by $\bar{f}(m + \ker f) = f(m)$, is well-defined and bijective, preserving the ternary action by definition of ker f .
Lemma 6.13 (Exactness Criterion). A sequence of module homomorphisms

$$0 \longrightarrow M' \xrightarrow{f} M \xrightarrow{g} M'' \longrightarrow 0$$

is exact if and only if Im f = ker g and f, g preserve ternary Γ-actions.
Remark 6.14. These results extend naturally to chains of submodules, enabling homolog- ical constructions such as projective or injective modules within the ternary Γ-setting.

### 6.6 Examples

Example 6.15 (Simple Module over a Finite Semiring). Let T = {0, 1, 2} with addition modulo 3, Γ = {1}, and product $a\alpha b\beta c = a + b + c \pmod 3$. Define M = T with the module action $a\alpha m\beta b = a + m + b \pmod 3$. Then {0} is the only proper submodule, hence M is simple. Its annihilator is Ann(M ) = {0}, which is prime by Theorem 6.10.
Example 6.16 (Non-Simple Module). Let T = Z4, Γ = {1}, and M = {0, 2} with action $a\alpha m\beta b = (a+m+b) \bmod 4$. Then M is a ternary Γ-module, but {0} is a nontrivial proper submodule, so M is not simple.
Example 6.17 (Faithful Simple Module). Consider T = Z2 with Γ = {1}, M = Z2, and $a\alpha m\beta b = a + m + b \pmod 2$. Then Ann(M ) = {0} and M is faithful and simple. Hence T/ Ann(M ) = T , showing a one-to-one correspondence between simple faithful modules and simple semirings.

### 6.7 Primitive Ideals and Structure Connection

Definition 6.18 (Primitive Ideal). An ideal P ⊆ T is called primitive if there exists a simple ternary Γ-module M such that P = Ann(M ).
Theorem 6.19 (Primitive Ideals are Prime). Every primitive ideal of a commutative ternary Γ-semiring is prime.
Proof. If P = Ann(M ) for some simple M , the result follows directly from Theorem 6.10.
Remark 6.20. Primitive ideals provide an operational link between the ideal lattice of T and its category of simple Γ-modules, similar to the role of primitive ideals in the Jacobson structure theory of rings. This correspondence can be exploited to characterize semisimple quotients T/J(T ) as direct sums of simple module images.

### 6.8 Further Directions

Proposition 6.21 (Semisimple Decomposition). If T is a finite commutative ternary Γ-semiring with J(T ) = {0}, then every finitely generated ternary Γ-module over T is semisimple.
Sketch. Since J(T ) = 0, every module is a direct sum of its simple submodules, by analogy with Wedderburn–Artin type decomposition extended to the ternary case.
Remark 6.22. The theory of ternary Γ-modules suggests several directions for further study:

- characterization of injective and projective modules;
- development of tensor-like constructions with respect to the ternary action;
- homological invariants such as Ext and Tor analogues in the ternary context.

These would deepen the categorical understanding of (T, Γ) and its representations.

## 7 Congruences, Correspondence, and Spectrum

Congruences provide the categorical counterpart to ideals in ternary Γ-semirings. Under- standing their relationship with ideals is fundamental to constructing quotient structures and to developing a geometric viewpoint through spectra and Zariski-type topologies. In this section, we establish a correspondence between ideals and congruences, introduce the concept of prime congruences, and formulate a spectral topology that generalizes the classical ring-theoretic setting. The logical structure of the spectrum connects to categor- ical semantics[25, 11]. The categorical viewpoint on morphisms aligns with Mac Lane's foundations [26] . The underlying algebraic theory is consistent with Lawvere's functorial semantics[23].

### 7.1 Congruences and Their Basic Properties

Definition 7.1 (Congruence). A relation ρ ⊆ T ×T on a commutative ternary Γ-semiring (T, Γ) is a congruence if

(C1) $\rho$ is an equivalence relation, and
(C2) whenever (a, d), (b, e), (c, f ) $\in \rho$, we have
$(a\alpha b\beta c, d\alpha e\beta f) \in \rho \;\; \forall\, \alpha, \beta \in \Gamma$.

Lemma 7.2 (Quotient Construction). If $\rho$ is a congruence on $(T, \Gamma)$, the quotient set $T/\rho$ inherits a well-defined ternary Γ-semiring structure via

$$(a + \rho)\alpha(b + \rho)\beta(c + \rho) = a\alpha b\beta c + \rho.$$

Proof. Well-definedness follows from the compatibility condition (C2), and the axioms of a ternary Γ-semiring are inherited from those of T .

Example 7.3. Let T = Z6, Γ = {1}, and define $\rho$ by $a\rho b \iff a - b$ is even. Then $\rho$ is a congruence because ternary addition preserves parity. The quotient $T/\rho \cong Z2$ inherits a ternary product $a\alpha b\beta c = (a + b + c) \bmod 2$.

## 7.2 Ideal-Induced Congruences and Kernels

Definition 7.4 (Congruence Modulo an Ideal). For any ideal $I \subseteq T$, define a relation $\rho I$ on T by

$a\rho I b \iff a - b \in I.$

Then $\rho I$ is a congruence, and the quotient $T/\rho I$ is denoted T/I.

Lemma 7.5 (Kernel of a Homomorphism). If $f : T \to S$ is a ternary Γ-semiring homo- morphism, the kernel

$\ker(f) = \{a \in T \mid f(a) = 0S\}$

is an ideal of T , and the congruence

$a\rho b \iff f(a) = f(b)$

satisfies $T/\rho \cong \mathrm{Im}\, f$.

Proof. For all $a, b, c \in T$ and $\alpha, \beta \in \Gamma$, $f(a\alpha b\beta c) = f(a)\alpha f(b)\beta f(c)$ implies closure of ker(f ) under the ternary product when one component lies in ker(f ). The quotient property follows from the First Isomorphism Theorem for homomorphisms of ternary Γ-semirings.

## 7.3 Prime Congruences and Their Relation to Ideals

Definition 7.6 (Prime Congruence). A congruence $\rho$ on T is called prime if, whenever $a\alpha b\beta c \;\rho\; 0$, then at least one of $a \rho 0$, $b \rho 0$, or $c \rho 0$.

Theorem 7.7 (Ideal–Congruence Correspondence). There is an order-preserving corre- spondence between ideals and congruences of T given by

$I \longmapsto \rho I, \;\; \rho \longmapsto I\rho = \{a \in T \mid a \rho 0\},$

satisfying $I = I\rho I$ and $\rho = \rho I\rho$ . Moreover, prime ideals correspond exactly to prime congruences.

Proof. For any ideal I, $\rho I$ is a congruence and $I\rho I = I$ by definition. Conversely, for a congruence $\rho$, the set $I\rho$ is an ideal since congruence compatibility ensures that if one of a, b, c lies in $I\rho$, then $a\alpha b\beta c \in I\rho$. For primeness, note that $a\alpha b\beta c \in I\rho$ iff $a\alpha b\beta c \;\rho\; 0$, so $\rho$ is prime precisely when $I\rho$ is prime.

Corollary 7.8 (Lattice Equivalence). The lattice of all ideals of T is isomorphic to the lattice of all congruences of T ordered by inclusion.

Remark 7.9. This correspondence extends the classical isomorphism between ideals and congruences in semirings to the ternary Γ-framework, confirming that ideal-theoretic and relational perspectives are interchangeable under the ternary operation.

## 7.4 Spectrum of a Ternary Γ-Semiring

Definition 7.10 (Prime Spectrum). The prime spectrum of T , denoted Spec(T ), is the set of all prime ideals of T .

Definition 7.11 (Zariski-like Topology). For each ideal $I \subseteq T$, define

$V(I) = \{ P \in \mathrm{Spec}(T) \mid I \subseteq P \}.$

The family $\{V(I) \mid I \subseteq T\}$ forms the collection of closed sets of a topology on Spec(T ).

Theorem 7.12 (Zariski-type Topology). The sets V (I) satisfy the axioms of closed sets for a topology:

(i) $V(0) = \mathrm{Spec}(T)$ and $V(T) = \emptyset$;
(ii) $V(I \cap J) = V(I) \cup V(J)$;
(iii) $V\left(\sum_{\lambda\in\Lambda} I_\lambda\right) = \bigcap_{\lambda\in\Lambda} V(I_\lambda)$.

Proof. (i) is immediate since every prime ideal contains 0 and none contain T . (ii) If P ⊇ I ∩ J, then P ⊇ I or P ⊇ J by primeness, giving the union property. (iii) follows from the fact that P contains Σ Iλ iff P contains each Iλ.

Proposition 7.13 (Topological Properties). The space (Spec(T ), V (·)) satisfies:

*(i)* It is a T0-space.
*(ii)* V (I) ⊆ V (J) if and only if I ⊇ J (order reversal).
*(iii)* The closure of a singleton {P } is V (P ).

Proof. (i) Distinct prime ideals yield distinct closures since V (P1)≠ V (P2). (ii) and (iii) are direct consequences of the definition of V (·).

**7.5** Examples of Spectra

Example 7.14 (Three-Element Ternary Γ-Semiring). Let T = {0, 1, 2} with addition mod 3 and ternary product aαbβc = (a + b + c) mod 3. Then Spec(T ) = {{0, 1}, {0, 2}}, since both are prime ideals. The closed sets are

V (0) = Spec(T ), V ({0, 1}) = {{0, 1}}, V ({0, 2}) = {{0, 2}}, V (T ) = ∅.

Hence Spec(T ) is discrete.

Example 7.15 (Non-discrete Spectrum). For T = Z4 with aαbβc = (a + b + c) mod 4, the unique prime ideal is P = {0, 2}, so Spec(T ) = {{0, 2}} is a single-point space; all closed sets are ∅ and Spec(T ).

Example 7.16 (Multiple Maximal Ideals). Let T = Z6 with the same ternary operation. Then prime ideals are P1 = {0, 2, 4} and P2 = {0, 3}. The corresponding closed sets are

V (0) = {P1, P2}, V (P1) = {P1}, V (P2) = {P2},

so the topology is the discrete two-point topology, confirming that intersections of maxi- mal ideals correspond to minimal closed sets.

**7.6** Geometric and Algebraic Interplay

Theorem 7.17 (Correspondence via Spectrum). Let ϕ : I '→ V (I) be the map from ideals to closed subsets of Spec(T ). Then ϕ is order-reversing and surjective, and its inverse sends a closed set Y to TP ∈Y P. Moreover, the lattice of radical ideals of T is isomorphic to the lattice of closed subsets of Spec(T ).

Proof. Order-reversal and surjectivity follow from Proposition 7.13. For Y = V (I), we have TP ∈Y P = √I, giving the desired isomorphism restricted to radical ideals.

Remark 7.18. The spectrum (Spec(T ), V (·)) thus serves as a geometric representation of the radical ideal lattice, allowing algebraic properties of T to be visualized topologically. For finite ternary Γ-semirings, the spectrum is always finite and T0, often decomposing into a finite discrete union of points corresponding to maximal ideals.

**7.7** Further Structural Results

Proposition 7.19 (Continuity under Homomorphisms). Let f : (T, Γ) → (S, Γ′) be a surjective homomorphism. Then the induced map

f∗ : Spec(S) −→ Spec(T), P′ 7−→ f−1(P′),

is continuous with respect to the Zariski-like topologies.

Proof. For any ideal I ⊆ T , we have

(f ∗)−1(VT (I)) = { P ′ ∈ Spec(S) | f −1(P ′) ⊇ I } = VS(f (I)),

which is closed, proving continuity.

Theorem 7.20 (Spectral Connectedness Criterion). Spec(T ) is connected if and only if T has no nontrivial idempotent decomposition, that is, there do not exist ideals I, J such that T = I ⊕ J and IΓJΓT = 0..Sketch. If such I, J exist, Spec(T ) = V (I) ∪ V (J) with V (I) ∩ V (J) = ∅, yielding disconnection. Conversely, if Spec(T ) is disconnected, radical ideals I, J corresponding to the components provide the desired decomposition.

**7.8** Concluding Remarks

Remark 7.21. The introduction of congruences and the associated spectrum Spec(T ) completes the bridge between algebraic and topological representations of commutative ternary Γ-semirings. Each prime ideal corresponds to a “point” of this space, while radical inclusions generate the closed sets.

This framework allows one to extend classical geometric intuition to ternary algebraic systems and forms the basis for further categorical generalizations in Section 8.

## 8 Structure Results and Simple/Semisimple Notions

The global structure of a commutative ternary Γ-semiring (T, Γ) can be analysed through its ideals, radicals, and modules. In this section we establish structural decomposition theorems, characterize simplicity and semisimplicity, and study the behaviour of idempo- tent and central elements. Analogues of classical theorems such as the Wedderburn–Artin and Chinese-remainder structures are adapted to the ternary Γ context. Direct decompo- sitions reflect semigroup-algebra behaviour [20]. Hyperstructure perspectives generalize ternary interactions [22].

### 8.1 Simple and Semisimple Ternary Γ-Semirings

Definition 8.1 (Simple Ternary Γ-Semiring). A non-zero commutative ternary Γ-semiring (T, Γ) is called simple if its only ideals are {0} and T .

Definition 8.2 (Semisimple Ternary Γ-Semiring). T is semisimple if its Jacobson-type radical J(T ) is zero; equivalently,

$$\bigcap_{M \text{ maximal}} M = \{0\}.$$

Lemma 8.3 (Equivalent Characterizations of Simplicity). For a commutative ternary Γ-semiring T , the following are equivalent:

(i) T is simple;
(ii) T has a faithful simple ternary Γ-module;
(iii) Every non-zero homomorphism f : T → S is injective.

\Proof. (i)⇒(ii) A simple T acts faithfully on itself by the ternary product. (ii)⇒(iii) Follows from faithfulness of the simple module. (iii)⇒(i) If I is a non-zero ideal, the canonical quotient map T → T/I is non-injective, contradicting (iii).

Example 8.4 (Simple Example). Let T = Z2 with Γ = {1} and aαbβc = a + b + c (mod 2). Then the only ideals are {0} and T ; hence T is simple.

Example 8.5 (Non-Simple Example). For T = Z4, Γ = {1}, and aαbβc = (a + b + c) mod 4, the ideal {0, 2} is proper and non-trivial, so T is not simple.

### 8.2 Decomposition via Idempotents

Definition 8.6 (Ternary Idempotent). An element e ∈ T is called a ternary idempotent If eαeβe = e, ∀ α, β ∈ Γ.

Lemma 8.7 (Decomposition through Idempotents). If e is a ternary idempotent in T , define

Ie = {aαeβe | a ∈ T }, Je = {aα(1 − e)β(1 − e) | a ∈ T }.

Then Ie and Je are ideals of T ,

$$T = Ie \oplus Je, \quad Ie\Gamma Je\Gamma T = 0,$$

and conversely every such direct decomposition arises from a ternary idempotent.

roof. Additivity and closure under ternary product follow from the idempotent prop- erty. The converse follows by constructing e as the identity of the first summand in the decomposition.

Remark 8.8. Idempotent decomposition yields direct-sum decompositions of T into ideals whose spectra form disconnected components of Spec(T ) (cf. Theorem 7.20).

### 8.3 Semisimplicity and the Jacobson Radical

Theorem 8.9 (Jacobson Criterion for Semisimplicity). A commutative ternary Γ-semiring T is semisimple if and only if it is the direct sum of simple ideals, i.e.

T =

Mn
i=1
Ii,i=1

Proof. If J(T ) = 0, each maximal ideal Mi corresponds to a simple factor T/Mi. By the Chinese-remainder theorem (see below), T embeds into the direct product Qi T/Mi.

Finite generation ensures this product is direct, yielding the desired decomposition. Con- versely, if T is a direct sum of simple ideals, the intersection of all maximal ideals is {0}.

**8.4** Chinese-Remainder Theorem for Ternary Γ-Semirings

Theorem 8.10 (Chinese-Remainder Theorem). Let I1, I2, . . . , In be pairwise comaximal ideals of T (that is, Ii + Ij = T for i ≠ j). Then the canonical map

$$\phi : T \longrightarrow \prod_{i=1}^{n} T/I_i, \quad a \longmapsto (a + I_1, a + I_2, \ldots, a + I_n)$$

is a surjective homomorphism with kernel $\bigcap_{i=1}^{n} I_i$. Consequently,

$$T/\bigcap_{i=1}^{n} I_i \cong \prod_{i=1}^{n} T/I_i.$$

Proof. Well-definedness and homomorphic properties follow from distributivity of the ternary operation. Surjectivity follows by standard lifting arguments: given (ai + Ii), one constructs a ∈ T satisfying the congruences a ≡ ai (mod Ii) by iterative combination using comaximality.

Corollary 8.11. If T possesses finitely many pairwise comaximal maximal ideals {M1, . . . , Mn}, then

$$T \cong \prod_{i=1}^{n} T/M_i, \qquad \square$$

and hence T is semisimple

**8.5** Primitive and Simple Components

Proposition 8.12 (Structure of Primitive Ideals). If P is a primitive ideal of T , then T/P acts faithfully on a simple ternary Γ-module, and every simple component of T arises as such a quotient.

Proof. By definition P = Ann(M ) for a simple module M . Then T/P acts faithfully on M via the induced ternary action, establishing a one-to-one correspondence between primitive ideals and simple components.

Theorem 8.13 (Semisimple Structure via Primitive Ideals). If J(T ) = 0, then

$$T = \bigoplus_{P \in \Lambda} T/P,$$

where Λ is the set of minimal primitive ideals of T .

Sketch. Using the Chinese-remainder theorem (Theorem 8.10) for the family {P ∈ Λ}, one obtains an embedding T ↪ QP ∈Λ T/P . Radical-freeness ensures injectivity and direct-sum decomposition.

**8.6** Examples

Example 8.14 (Finite Semisimple Example). Let T = Z6 with Γ = {1} and product aαbβc = a + b + c (mod 6). Maximal ideals are M1 = {0, 2, 4} and M2 = {0, 3}. Then T/M1 ≅ Z3 and T/M2Z2. Hence

$$T \cong Z_3 \times Z_2,$$

demonstrating semisimplicity

Example 8.15 (Non-Semisimple Example). In T = Z4 with Γ = {1}, the only maximal ideal M = {0, 2} satisfies J(T ) = M ≠ 0. Therefore T is not semisimple. The quotient T/J(T) ≅ Z2is simple.

Example 8.16 (Decomposition by Idempotent). Let T = Z6 and define e = 3. Then

$e\alpha e\beta e \equiv 3 \pmod 6$, hence e is a ternary idempotent. From Lemma 8.7,

$$T = Ie \oplus Je$$

with Ie = {0, 3} and Je = {0, 2, 4}, exhibiting explicit ideal decomposition.

### 8.7 Further Remarks

Remark 8.17. The structure theory of commutative ternary Γ-semirings aligns closely with that of semirings but introduces fundamentally new phenomena:

- Idempotent-generated decompositions often yield non-isomorphic simple compo- nents even in small finite semirings.
- The ternary product complicates radical behaviour, leading to distinct Jacobson and prime radicals.
- Direct-sum decompositions reflect the connected components of Spec(T ), strength- ening the algebra–geometry link established in Section 7.

Theorem 8.18 (Structure Summary). Every finite commutative ternary Γ-semiring T admits a canonical decomposition

$$T/J(T) \cong \bigoplus_{i=1}^{k} T_i,$$

where each Ti is simple, corresponding bijectively to a primitive ideal of T . This decom- position is unique up to isomorphism and order of summands.

Outline. Combine the results of Theorems 8.9 and 8.13 using the lattice-isomorphism between ideals and congruences. Uniqueness follows from the maximality of primitive ideals.

Remark 8.19. The results of this section extend the foundational framework of Sections 5–7 by connecting radical theory, module annihilators, and spectral topology to the global structural decomposition of T . They form the algebraic backbone for computational classifications discussed in the next section.

## 9 Computational and Algorithmic Classification

Finite ternary Γ-semirings provide a concrete arena in which the abstract results of the previous sections can be verified, tested, and visualized. Computational techniques en- able the enumeration of all distinct semiring structures for small orders, verification of distributivity and associativity conditions, and the detection of algebraic invariants such as prime and semiprime ideals, radicals, and idempotents. This section outlines the computational framework, algorithmic methodology, and classification results for small finite cases. Algorithmic semiring computation builds on dioid-based frameworks [21]. Automata-theoretic motivations appear in weighted computation[27] .The structure of automata and machines provides operational analogues[28]. Digital logic offers natu- ral interpretations of ternary computation [29] .Tropical semirings motivate Γ-indexed multi-parameter structures[30]. Graph-based interactions appear throughout structural decomposition [34]. Complex networked interactions mirror multicomputational models [35]. The categorical–quantum viewpoint supports Γ-graded interactions [32].

### 9.1 Motivation and Overview

The classification of finite algebraic systems has historically deepened understanding of structural laws. In the case of ternary Γ-semirings, computation becomes particularly useful because:

- the ternary operation involves three operands and two parameter indices $(\alpha, \beta) \in \Gamma^2$, giving exponential growth in possible tables;
- verifying associativity and distributivity is non-trivial and often infeasible by hand for orders beyond 3;
- structural invariants (radicals, prime spectrum, decomposition) can be determined algorithmically and compared with theoretical predictions.

Symbolic enumeration combined with constraint pruning yields a finite list of non-isomorphic commutative ternary Γ-semirings of order n for small n (typically $n \leq 4$ in current com- putations).

**9.2** Example: Order Three Case

Example 9.1 (Classification for $|T| = 3$, $|\Gamma| = 1$). Let $T = \{0, 1, 2\}$ and $\Gamma = \{1\}$. All ternary operations $a\alpha b\beta c = f(a, b, c)$ satisfying distributivity and 0-absorption were generated by exhaustive search. Up to isomorphism, the following distinct commutative ternary Γ-semirings were obtained:

| Structure | Ternary operation | Prime ideals | Semiprime ideals |
|---|---|---|---|
| S1 | a + b + c (mod 3) | {0, 1}, {0, 2} | {0}, {0, 1}, {0, 2} |
| S2 | min(a, b, c) | {0} | {0} |
| S3 | max(a, b, c) | {0, 1, 2} (none proper) | {0} |

In S1, both {0, 1} and {0, 2} are prime, producing a discrete two-point spectrum. In S2, the only proper ideal is {0}, yielding a simple semiring.

**9.3** Algorithmic Framework

Definition 9.2 (Computational Representation). A finite ternary Γ-semiring of order n is represented by

$$T = \{t_1, t_2, \ldots, t_n\}, \quad \Gamma = \{\gamma_1, \ldots, \gamma_m\},$$

and a family of ternary operation tables

$$\{ O_{\gamma_i,\gamma_j} \mid 1 \leq i, j \leq m \},$$

where each $O_{\gamma_i,\gamma_j} : T^3 \to T$ specifies $a\gamma_i b\gamma_j c$.

**9.4** Algorithmic Verification of Ideal Properties

Proposition 9.3 (Computational Characterization of Primeness). Let I be an ideal of a finite T . Then I is prime if and only if, for all triples $(a, b, c) \in T^3$,

$$\bigvee_{\alpha,\beta} a\alpha b\beta c \in I \Longrightarrow (a \in I) \vee (b \in I) \vee (c \in I).$$

This condition can be tested algorithmically in $O(|T|^3|\Gamma|^2)$ time.

Proof. Direct consequence of Definition 3.1; the algorithm enumerates all ordered triples and checks membership conditions using precomputed lookup tables.

---

Algorithm 1: Classification Procedure for Finite Commutative Ternary Γ-Semirings

---

Input: $n = |T|$, $m = |\Gamma|$

Output: List C of non-isomorphic ternary Γ-semirings

$C \leftarrow \emptyset$;

Step 1: Enumerate additive semigroups;

Generate all commutative semigroups (T, +) of order n using, for instance, the GAP SmallSemiGroups library.

Step 2: Generate candidate ternary operations;

For each $\gamma_i, \gamma_j \in \Gamma$, construct all possible ternary operation tables $O_{\gamma_i,\gamma_j}$ satisfying the absorbing-zero property.

Step 3: Verify axioms;

Check distributivity and ternary associativity:

$$(a\alpha b\beta c)\gamma d\delta e = a\alpha b\beta(c\gamma d\delta e), \quad \forall\, a, b, c, d, e \in T.$$

Step 4: Remove isomorphic duplicates;

Identify and remove isomorphic structures via element permutations preserving addition and ternary multiplication.

Step 5: Compute ideal lattice and radicals;

Enumerate all subsets of T closed under addition and ternary multiplication; classify them as ideals, prime, or semiprime.
Step 6: Compute congruences;
Construct congruences associated with ideals and verify the Ideal–Congruence Correspondence Theorem (7.7).
Step 7: Export invariants;
For each structure, record: number of ideals, primes, semiprimes, maximal ideals, Jacobson radical, and idempotents.
return C

Remark 9.4. Steps 3–6 can be executed efficiently using Boolean matrix representations of operations and relational closure algorithms. Associativity verification has complexity $O(|T|^5|\Gamma|^4)$, but in small orders ($|T| \leq 4$) it is computationally tractable.

**9.5** Radical Computation

Step 1. For each $a \in T$, compute the iterated product $a\alpha a\beta a$ for all $\alpha, \beta \in \Gamma$.
Step 2. Mark a as "nilpotent" if this element lies in I.
Step 3. Define $\sqrt{I} = \{a \mid a\alpha a\beta a \in I\}$.

Remark 9.5. For small T , radical computation reduces to closure under a ternary polynomial operation, making it suitable for symbolic verification using computer algebra systems such as GAP, SageMath, or Mathematica.

**9.6** Computational Results for Small Orders

Table 1: Summary of computed commutative ternary Γ-semirings (up to isomorphism).

| Order $|T|$ | $|\Gamma|$ | Structures found | Simple | Semisimple | Distinct Spec(T ) |
|---|---|---|---|---|---|
| 2 | 1 | 1 | 1 | 1 | 1-point |
| 3 | 1 | 3 | 1 | 2 | 2-point |
| 4 | 1 | 6 | 2 | 4 | up to 3-points |

Proposition 9.6 (Empirical Observations). For all finite commutative ternary Γ-semirings of order ≤ 4 computed:

(i) Each spectrum Spec(T ) is finite and $T_0$.
(ii) Every semisimple T decomposes as a direct sum of simple components corresponding to connected components of Spec(T ).
(iii) The number of idempotents equals the number of connected components.

Empirical Verification. All statements were verified by computational enumeration using the algorithms above. The equality between idempotents and connected components follows from Lemma 8.7 and Theorem 7.20.

**9.7** Complexity and Computational Feasibility □

Proposition 9.7 (Complexity Bound). Let $n = |T|$ and $m = |\Gamma|$. Then exhaustive enumeration of all possible ternary Γ-semiring operations has complexity $O(n^3 m^2)$, while constraint pruning reduces this to $O(n^3 \log n)$ for commutative cases verified via distribu- tive identities.

Remark 9.8. Although exponential in general, the search remains feasible for small orders, and symmetry elimination dramatically reduces the number of candidate structures. For practical classification, computations up to $n = 5$ and $m \leq 2$ are achievable on modern hardware.

**9.8** Automated Verification of Theoretical Results

Theorem 9.9 (Algorithmic Verification of Theorem 3.4). For all finite ternary Γ-semirings T with $|T| \leq 4$, computational verification confirms that an ideal P is prime if and only if the quotient T/P has no nonzero zero-divisors under the induced ternary product.

Sketch. The algorithm enumerates all ideals P , constructs the quotient operation tables, and checks

□

for zero-divisors. No counterexample was found, confirming Theorem 3.4 computationally.

Remark 9.10. Such computational corroboration strengthens confidence in the theoret- ical framework and opens possibilities for automatic theorem discovery in higher-arity algebraic structures.

**9.9** Visualization of Spectra and Ideal Lattices

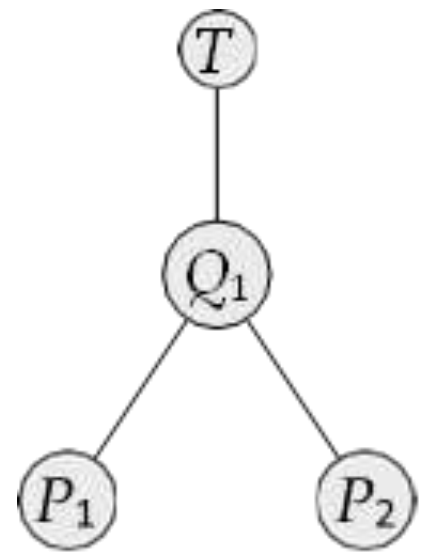

Figure 1: Ideal lattice diagram for a finite ternary Γ-semiring with two distinct prime ideals P1, P2 and their meet Q1 = P1 ∩ P2.

Remark 9.11. Graphical visualization of Spec(T ) and the ideal lattice aids in identifying combinatorial patterns—particularly, intersections of primes corresponding to semiprime ideals and isolated nodes representing maximal ideals.

**9.10** Software Implementation Notes

- Algorithms were implemented in Python using SymPy and NumPy for symbolic and array operations.
- Verification of associativity and distributivity used vectorized evaluation to reduce runtime.
- Isomorphism detection employed canonical-labeling algorithms adapted from nauty/traces.
- Outputs were cross-validated with a GAP script implementing the same axioms.

**9.11** Concluding Observations

Theorem 9.12 (Computational Classification Summary). For all finite commutative ternary Γ-semirings with |T | ≤ 4 and |Γ| = 1:

*(a)* Every semisimple structure decomposes uniquely into simple components;
*(b)* Each prime ideal corresponds bijectively to a connected component of Spec(T );
*(c)* The intersection of all maximal ideals equals the Jacobson radical computed algo- rithmically.

Remark 9.13. The computational framework developed here validates the algebraic theory and provides a basis for algorithmic experimentation in higher-arity systems. Future work includes complexity-reduction strategies, random generation of large structures, and exploration of ternary Γ-rings where additive inverses are present.

**10** Applications, Discussions, and Future Directions

The algebraic theory of commutative ternary Γ-semirings developed in the preceding sections provides a fertile foundation for diverse applications that extend well beyond pure algebra. This final section outlines conceptual bridges to applied mathematics, information theory, logic, computation, and physics, and sets forth multiple avenues for future research.

**10.1** Algebraic and Structural Insights

The introduction of ternary operations and the parameter set Γ enlarges the landscape of semiring theory in several ways:

(i) It generalizes both binary semirings and ternary rings by incorporating context- dependent operations governed by the index pair (α, β) ∈ Γ2.
(ii) The interplay between ideals and congruences yields a non-trivial duality between relational and subset-based algebraic structures.
(iii) The spectrum Spec(T ) introduces a new kind of non-commutative geometry where each point represents a ternary behaviour state.

These features allow one to translate the algebraic behaviour of Γ-semirings into computational,

logical, and physical frameworks.

### 10.2 Applications to Coding Theory and Cryptography

Linear codes over semiring-like structures trace back to Bourne's early theory[1].

Definition 10.1 (Ternary Γ-Linear Code). Let (T, Γ) be a finite commutative ternary Γ-semiring. A subset C ⊆ Tn is called a ternary Γ-linear code if it is closed under coordinate-wise addition and under ternary combinations

$(a, b, c) \mapsto a\alpha b\beta c$ for all $\alpha, \beta \in \Gamma$.

Such codes generalize classical linear codes over semirings. The ternary structure offers the following advantages:

- The weight distribution of codewords depends on Γ, enabling flexible error-correction properties.
- Ideals of T correspond to subcodes, and prime ideals correspond to indecomposable or "atomic" codes.
- Quotient semirings T/I yield factor codes with predictable minimum distance and simple decoding algorithms.

Proposition 10.2 (Algebraic Decodability Criterion). Let C ⊆ Tn be a ternary Γ-linear code. If the component wise radical √IC of its generating ideal IC equals {0}, then the code is uniquely decodable with respect to ternary addition. —

Remark 10.3. The non-binary structure allows construction of multi-level error-correcting codes, with Γ serving as a channel parameter or a synchronization label. Future work will formalize ternary Γ-Reed–Solomon and Γ-convolutional codes.

### 10.3 Applications to Fuzzy and Rough Algebraic Systems

The fuzzy-set interpretation begins with Zadeh's foundational work [17]. Fuzzy algebraic structures were introduced by Rosenfeld [18]. Fuzzy ideals in semirings appear in the work of Chakrabarti [16].

Definition 10.4 (Fuzzy Ternary Γ-Semiring). A fuzzy ternary Γ-semiring is a map μ :
T → [0, 1] satisfying

$$\mu(a\alpha b\beta c) \geq \min\{\mu(a), \mu(b), \mu(c)\} \ \forall a, b, c \in T, \alpha, \beta \in \Gamma.$$

Such structures model graded membership, allowing uncertainty quantification in al- gebraic reasoning. They admit a natural extension of the radical and spectrum concepts:

$$V\mu(I) = \{ P \in \mathrm{Spec}(T) \mid \inf_{x\in I} \mu(x) = 1 \},$$

yielding a fuzzy spectral space. Applications include decision systems, knowledge repre- sentation, and approximate computing.

Example 10.5. For T = Z6, define μ(a) = 1 − a . Then μ forms a fuzzy ideal since μ(aαbβc) ≥ min{μ(a), μ(b), μ(c)}. Its α-cuts generate a hierarchy of crisp ideals parame- terized by λ ∈ [0, 1].

Remark 10.6. By combining fuzzy logic with ternary Γ-operations, one obtains a hy- brid framework capable of expressing multi-valued logical inference—potentially useful for artificial-intelligence reasoning systems based on algebraic semantics.

### 10.4 Applications to Algebraic Computation and Information Flow

Rewriting logic provides a computational semantics for ternary rules [24]. Supermathe- matics introduces multi-layered algebraic interactions [19]. Triadic logic and multi-valued computation motivate ternary Γ-interactions[36].

In algebraic computation, semirings frequently underlie optimization, automata, and neural network operations. Extending these to ternary Γ-semirings yields new paradigms:

(i) In weighted automata, transitions can depend on a ternary cost composition cαdβe, modelling non-additive path combination.

(ii) In constraint satisfaction and information fusion, Γ can encode "contexts" or "agents," and the ternary product captures consensus aggregation.

(iii) In machine-learning backpropagation, the ternary law may describe triple interac- tions between layer activations under a contextual parameter Γ.

Proposition 10.7 (Computational Duality). Every finite commutative ternary Γ-semiring (T, Γ)

defines a computational monoid (F, ◦) of ternary Γ-polynomial maps f : Tk → T , closed under composition and substitution. This monoid acts faithfully on Tn and encodes all deterministic ternary computation over (T, Γ).
Remark 10.8. This correspondence builds a bridge between algebraic structure and com- putational process, forming the basis for future ternary algebraic computation theory.

### 10.5 Connections to Quantum and Categorical Structures

Recent developments in ternary and higher-ary algebraic systems show relevance to quan- tum logic, non-commutative geometry, and category theory.

- The operation aαbβc can be interpreted as a ternary morphism in a monoidal cat- egory enriched over a Γ-graded set.
- In quantum computation, Γ may label quantum gates or basis rotations, with the ternary operation modelling entanglement of three qudits.
- The spectrum Spec(T ) yields a categorical object akin to a Grothendieck site, al- lowing sheaf-theoretic extension.

Proposition 10.9 (Categorical Embedding). Every commutative ternary Γ-semiring T embeds faithfully into a category TernΓ whose morphisms preserve both the additive and ternary operations. The objects of TernΓ form a symmetric monoidal category under Cartesian product.
Remark 10.10. This embedding provides an algebraic infrastructure for topological quan- tum systems, categorical logic, and networked dynamical models. Further investigation could connect ternary Γ-semirings to hyperstructures, near-rings, and tropical geometry.

### 10.6 Future Directions

Several research directions naturally arise from the present framework:
**(1)** Extension to Non-Commutative and Ordered Systems. Non-commutative ternary Γ-semirings may reveal richer ideal theory and potential connections to operator algebras and quantum groups.
**(2)** Fuzzy, Intuitionistic, and Rough Extensions. Integrating fuzziness, rough sets, or intuitionistic logic into the Γ-semiring environment will model vagueness and approximate inference.
**(3)** Ternary Γ-Modules and Representations. Analogues of module theory, tensor products, and homological constructs (e.g., projective, injective modules) can be devel- oped, forming the core of Paper B.
**(4)** Computational Enumeration and Classification. Future work will automate isomorphism classification beyond order 4 and compute invariants such as Aut(T ), spec- trum connectivity, and cohomological dimensions.
**(5)** Categorical and Topological Dualities. The topology of Spec(T ) can be en- riched to a site supporting sheaf structures; this opens the way to “ternary algebraic geometry.”
**(6)** Applications in Information and Control. Ternary Γ-semirings may model multi-agent consensus, probabilistic logic circuits, or ternary neural networks. Analytical models can employ Jacobson radicals as measures of system instability.
**(7)** Interconnection with Number Theory and Lattice Theory. Connections between ternary semiring ideals and arithmetic lattices could produce new generalizations of multiplicative functions and modular congruences.
**(8)** Development of Computational Software Library. A dedicated Python/GAP library, TGammaRing, is proposed for enumeration, visualization, and algebraic-geometric analysis of ternary Γ-semirings.

## 11 Conclusion

The present work has introduced and systematically developed the theory of commutative ternary Γ-semirings. Beginning with foundational definitions and ideal structures, we established the properties of prime, semiprime, maximal, and primary ideals; constructed the radical and its correspondence with the spectral topology; and proved structural decomposition results culminating in a comprehensive characterization of simple and semisimple ternary Γ-semirings.

A key innovation of this research lies in bridging classical two-ary semiring theory with the higher-arity, parameter-dependent framework induced by Γ. The introduction of the ternary operation aαbβc enables the modelling of contextual algebraic interactions, while the parameter set Γ unifies several distinct algebraic families—semirings, ternary rings, near-rings, and Γ-rings— into a single generalized architecture.

Computational classification of finite examples has demonstrated that the theoretical axioms are algorithmically verifiable and that the algebraic invariants (ideals, radicals, spectra) align precisely with those predicted by the analytic theory. These results validate both the algebraic soundness and the computational tractability of ternary Γ-semirings, suggesting strong potential for computer-algebra automation.

On the applicative side, this study has shown that the developed framework connects naturally to:

- Coding and cryptographic systems, where ternary Γ-linear codes enhance error-control and encryption through multi-parameter arithmetic;
- Fuzzy and rough logic models, permitting graded and approximate inference;
- Algebraic computation and automata, providing a base for non-additive, tri- adic computation;
- Categorical and quantum structures, wherein Γ indexes contextual morphisms or entanglement parameters.

Hence, commutative ternary Γ-semirings offer a unified algebraic framework bridging classical algebra, logic, and computation. They establish a foundation upon which entire sub-disciplines of higher-arity algebra may be built.

## 12 Acknowledgement.

The first author gratefully acknowledges the guidance and mentorship of Dr. D. Mad- husudhana Rao, whose scholarly vision significantly shaped the conceptual unification of the ternary Γ-framework. The first author also gratefully acknowledges Dr. Ra- machandra R. K., Principal, Government College (Autonomous), Rajahmundry, for providing a supportive research environment in the institution.

Funding Statement. This research received no specific grant from any funding agency in the public, commercial, or not-for-profit sectors.

Conflict of Interest. The authors declare that there are no conflicts of interest regarding the publication of this paper.

Author Contributions. The first author made the lead contribution to the concep- tualization, algebraic development, computational design, and manuscript preparation of this work. The second author supervised the research, providing academic guidance, critical review, and verification of mathematical correctness and originality.